\documentclass[11pt]{article}
\usepackage{cite}
\usepackage{mathrsfs}
\usepackage{amsfonts}
\usepackage{amsmath}
\usepackage{amsfonts,amssymb,color}
\usepackage{dsfont}
\usepackage{curves}
\usepackage{mathrsfs}
\usepackage{pifont}
\usepackage{amssymb}
\usepackage{latexsym,amsmath,amssymb,amsfonts,epsfig,graphicx,cite,psfrag}
\usepackage{eepic,color,colordvi,amscd}

\newtheorem{theorem}{Theorem}[section]

\newtheorem{lemma}{Lemma}[section]
\newtheorem{corollary}{Corollary}[section]

\newtheorem{claim}{Claim}[section]

\newcommand{\qed}{\hfill\rule{0.5em}{0.809em}}

\def\emptyset{\mbox{{\rm \O}}}
\def\bar{\overline}

\renewcommand{\baselinestretch}{1.15}
\def\qed{\hfill \rule{4pt}{7pt}}

\def\pf{\noindent {\it Proof. }}

\begin{document}
	
	\title{Structure and coloring of a family of  ($P_7, C_5$)-free graphs}
	\author{Ran Chen\footnote{Email: 1918549795@qq.com},  \; Baogang  Xu\footnote{Email: baogxu@njnu.edu.cn. Supported by NSFC 11931006}\\\\
		\small Institute of Mathematics, School of Mathematical Sciences\\
		\small Nanjing Normal University, 1 Wenyuan Road,  Nanjing, 210023,  China}
	%\small $^2$School of Mathematics, Southeast University, 2 SEU Road, Nanjing, 211189, China}
\date{}

\maketitle

\vskip -50pt

\begin{abstract}
Let $P_t$ and $C_t$ be a path and a cycle on $t$ vertices, respectively. In 2021, Choudum {\em et al.} [Disc. Math. 344 (2021) 112244] determined the structures of  $(P_7,C_7,C_4$, diamond)-free and $(P_7,C_7,C_4$, gem)-frees, and gave correspondingly tight upper bounds to the chromatic numbers of these graphs.  In this paper, we study the structure of $(P_7, C_5$, kite, paraglider)-free graphs, which is a superfamily of $(P_7, C_5$, diamond)-free graphs. We show that there is a unique connected imperfect $(P_7, C_5$, kite, paraglider)-free graph with $\delta(G)\geq\omega(G)+1$, which has no  clique cutsets, no universal cliques, and no  pair of vertices of which one's neighborhoods contains the other's. As a consequence, we show that $(P_7, C_5$, kite, paraglider)-free graphs are $\chi$-polydet with a binding function $\omega(G)+1$. Where a {\em diamond} (resp. {\em gem}) consists of a $P_3$ (resp. $P_4$) and a new vertex adjacent to all vertices of the $P_3$ (resp. $P_4$), a {\em kite} consists of a  $P_4$ and a new vertex adjacent to consecutive three vertices of the $P_4$, and a {\em paraglider} consists of a  $C_4$ and a new vertex adjacent to three vertices of the $C_4$.
	
	\begin{flushleft}
	\end{flushleft}
	
\end{abstract}

\section{Introduction}\label{1}
We only consider finite and simple graphs, and follow \cite{BM08} for undefined notations. As usual, $P_t$ and $C_t$ denote a path and a cycle on $t$ vertices, respectively. Let $G$ be a graph, and $X$ and $Y$ be two disjoint subsets of $V(G)$. Let $E(X, Y)$ denote the set of edges between $X$ and $Y$. We say that $X$ is complete (resp. anticomplete) to $Y$ if $|E(X, Y)|=|X||Y|$ (resp. $|E(X, Y)|=0$).

For $u$, $v\in V(G)$, we write $u\sim v$ if $uv\in E(G)$, and write $u\not\sim v$ if $uv\not\in E(G)$.  Let $N_G(v)=\{x\;|\; x\in V(G), x\sim v\}$, $d_G(v)=|N_G(v)|$. Let $N_G(X)=(\cup_{x\in X} N_G(x))\setminus X$, and let $G[X]$ be the subgraph of $G$ induced by $X$. If it does not cause any confusion, we usually omit the subscript $G$.  Let $\delta(G)$ denote the minimum degree of $G$. For $X\subset V(G)$ and $v\in V(G)$, let $N_X(v)=N(v)\cap X$.

Let $G$ and $H$ be graphs. We say that $G$ contains $H$ if $H$ is isomorphic to an induced subgraph of $G$, and say that $G$ is $H$-{\em free} if it does not contain $H$.
For a family $\{H_1,H_2,\cdots\}$ of graphs, $G$ is $(H_1, H_2,\cdots)$-free if $G$ is $H$-free for each $H\in \{H_1,H_2,\cdots\}$.

A {\em clique} (resp. {\em stable set}) of $G$ is a set of mutually adjacent (resp. non-adjacent) vertices in $G$. The {\em clique number} $\omega(G)$ of $G$ is the maximum size of a clique in $G$. For a positive integer $k$, a $k$-{\em coloring } of $G$ is a function $\phi: V(G)\rightarrow \{1,\cdots,k\}$, such that $\phi(u)\ne \phi(v)$ if $u\sim v$. The {\em chromatic number} $\chi(G)$ of $G$ is the minimum integer $k$ for which $G$ has a $k$-coloring. A graph is {\em perfect} if all its induced subgraphs $H$ satisfy $\chi(H)=\omega(H)$. An induced cycle of length at least 4 is called  a {\em hole}, and its complement is called an {\em antihole}. A hole or antihole is {\em odd} or {\em even} if it has odd or even number of vertices. A {\em k-hole} is a hole of length $k$. Here below is the famous {\em Strong Perfect Graph Theorem} \cite{CRST2006}.
\begin{theorem}\label{perfect}{\em\cite{CRST2006}}
	A graph is perfect if and only if it is (odd hole, odd antihole)-free.
\end{theorem}

Let ${\cal G}$ be a family of graphs.  If there exists a function $f$ such that $\chi(G)\leq f(\omega(G))$ for every graph $G\in {\cal G}$, then we say that ${\cal G}$ is $\chi$-{\em bounded}, and call $f$ a {\em binding function} of ${\cal G}$ \cite{G75}.

It is known that $P_4$-free graphs are  perfect \cite{S74}. Improving the upper bound $(t-1)^{\omega(G)-1}$ of Gy\'{a}rf\'{a}s \cite{G75}, Gravier {\em et al} \cite{GHM03} proved that $\chi(G)\leq (t-2)^{\omega(G)-1}$ for all  $P_t$-free graphs. The question that whether there exists a polynomial binding functions to $P_t$-free graphs \cite{S16} is still open.  The best binding function for $P_5$-free graphs is due to Scott {\em et al} \cite{SSS2023}, who proved that every $P_5$-free graph $G$ satisfies $\chi(G)\leq\omega(G)^{\log_2(\omega(G))}$ for $\omega(G)\ge3$. Interested readers are referred to \cite{RS2004, SR2019, SS18} for more results and problems on $\chi$-bounded problems, and to \cite{G22} for most recent progresses on binding functions of subclasses of $P_5$-free graphs.

%Choudum and Karthick \cite{CKS07} conjectured that there exists a constant $c$ such that $\chi(G)\leq c\omega(G)^2$ if $G$ is $P_5$-free.

\begin{figure}[htbp]
	\begin{center}
		\includegraphics[width=10cm]{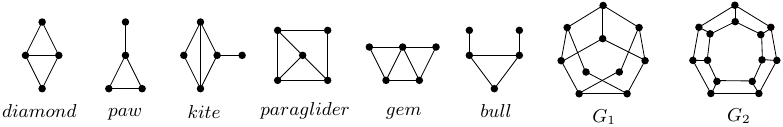}
	\end{center}
	\vskip -15pt
	\caption{\small Diamond, paw, kite, paraglider, gem, bull, $G_1$ and $G_2$.}
	\label{fig-1}
\end{figure}

There are also quite a lot of results on the binding functions of subclasses of $P_t$-free graphs when $t\ge 6$. We define more configurations such as {\em diamond}, {\em paw}, {\em kite}, {\em paraglider}, {\em gem}, and {\em bull}, as shown in Figure~\ref{fig-1}.

Choudum {\em et al.} \cite{CKS07} proved that  $\chi(G)\leq 8\omega(G)$ if $G$ is ($P_6$, gem)-free. Cameron {\em et al.} \cite{CHM18} showed that $\chi(G)\leq \omega(G)+3$ if $G$ is ($P_6$, diamond)-free, and Goedgebeur {\em et al.}  \cite{GHJ2021} improved the upper bound to $\max\{6,\omega(G)\}$. Improving an upper bound $\lfloor\frac{3}{2}\omega(G)\rfloor$ of Gasper and Huang \cite{GSH17}, Karthick and Maffray \cite{KM19} proved that $\chi(G)\leq\lceil\frac{5}{4}\omega(G)\rceil$ for $(P_6, C_4)$-free graph $G$. Mishra \cite{M21} proved that every ($P_t$, diamond, bull)-free graph $G$ has $\chi(G)\leq \max\{2t-4,\omega(G)\}$, and every ($P_7$, diamond, bull)-free graph $G$ satisfies $\chi(G)\leq \max\{7, \omega(G)\}$. On $(P_7, C_4, C_5)$-free graphs $G$, Huang \cite{H22} improved an upper bound $\frac{3}{2}\omega(G)$ of  Cameron {\em et al.} \cite{CHPS20}, and proved that $\chi(G)\leq\lceil\frac{11}{9}\omega(G)\rceil$.

Choudum {\em et al.} \cite{CKB21} proved that $\chi(G)\leq \max\{3,\omega(G)\}$ if $G$ is $(P_7,C_7,C_4$, diamond)-free, and $\chi(G)\leq 2\omega(G)-1$ if $G$ is $(P_7,C_7,C_4$, gem)-free. Chen {\em et al.} \cite{CWX2023} extended both upper bounds to larger families of graphs by dropping the requirement of $C_7$-free. Chudnovsky and Stacho \cite{MCJS2018} studied  $P_8$-free graphs, and proved that $\chi(G)\le 3$ if $G$ is $(P_8,C_3, C_i)$-free with $i=4$ or 5, and $\chi(G)\le 3$ if $G$ is $(P_8,C_4, C_5)$-free unless $G$ has subgraphs isomorphic to one of five well-defined small graphs.

It is known that $\chi(G)\le \omega(G)+1$ if $G$ is ($P_5$, diamond)-free \cite{BR1998}, and $\chi(G)\leq \omega(G)+3$ if $G$ is  ($P_6$, diamond)-free \cite{CHM18}. A still open question of Schiermeyer \cite{SR2023} asked that whether there exists a constant $c$ such that $\chi(G)\le \omega(G)+c$ for every $(P_7$, diamond)-free graph.
\begin{figure}[htbp]
	\begin{center}
		\includegraphics[width=2cm]{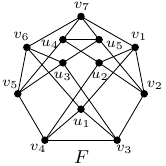}
	\end{center}
	\vskip -15pt
	\caption{\small The graph $F$.}
	\label{fig-2}
\end{figure}

Notice that the class of diamond-free graphs is a subclass of (paraglider, kite)-free graphs. In this paper, we study the structure of imperfect ($P_7, C_5,$ paraglider, kite)-free graphs, and prove the following Theorem~\ref{kite, paraglider}. A {\em clique cutset} of $G$ is a clique $K$ in $G$ such that $G-K$ has more components than $G$. Let $u$ and $v$ be two non-adjacent vertices. We call $(u, v)$ a {\em comparable pair} if $N(u)\subseteq N(v)$ or $N(v)\subseteq N(u)$. A clique $K$ is called a {\em universal clique} if $K$ is complete to $V(G)\setminus K$. Let $F$ be the graph as shown in Figure~\ref{fig-2}.
\begin{theorem}\label{kite, paraglider}
Let $G$ be a connected imperfect $(P_7, C_5$, kite, paraglider)-free graph with $\delta(G)\geq\omega(G)+1$. If $G$ is not isomorphic to $F$, then $G$ has a clique cutsets, or  a comparable pair, or a universal clique.
\end{theorem}

The following four classes of graphs show that the requirements $P_7$-free, $C_5$-free, kite-free, and paraglider-free are all necessary in Theorem~\ref{kite, paraglider}. Let $G_1$ and $G_2$ be the graphs as shown in Figure~\ref{fig-1}. Let $G_3$ be a graph whose vertex set is partitioned into seven cliques $Q_1$, $\cdots$, $Q_7$ such that for each $i$ mod 7, $|Q_i|=|Q_{i+1}|\geq2$, and $Q_i$ is complete to $Q_{i+1}\cup Q_{i-1}$ and anticomplete to $Q_{i+2}\cup Q_{i+3}\cup Q_{i-2}\cup Q_{i-3}$. Let $G_4$ be the complement of $C_7$. One can easily check that $G_1$ is $(P_7$, kite, paraglider)-free, $G_2$ is $(C_5$, kite, paraglider)-free, $G_3$ is $(P_7, C_5$, paraglider)-free, and $G_4$ is $(P_7,C_5$, kite)-free. But none of $G_1$, $G_2$, $G_3$ and $G_4$ satisfies Theorem~\ref{kite, paraglider}.

It is easy to check that a subset $S$ of $V(G)$ is a universal clique if and only if each vertex of $S$ has degree $|V(G)|-1$. If a graph $G$ has a comparable pair $(x, y)$ with $N(x)\subseteq N(y)$, then $\chi(G)=\chi(G-x)$. Notice that the graph $F$ is 3-colorable, and clique cutsets and universal cliques are reducible in coloring of graphs. Since clique cutsets can be found polynomially \cite{RET1985}, and since perfect graphs can be polynomially colored \cite{GLS1988}, as immediate consequence of Theorems~\ref{kite, paraglider}, we have the following corollary immediately by a simple induction on $|V(G)|$, which means the class of ($P_7, C_5$, kite, paraglider)-free graphs is $\chi$-polydet. Where a class ${\cal G}$ of graphs is said to be $\chi$-polydet \cite{SR2019} if ${\cal G}$ has a polynomial binding function $f$ and there exists a polynomial time algorithm to determine a coloring of $G\in {\cal G}$ with at most $f(\omega(G))$ colors.
\begin{corollary}\label{kite, paraglider*}
Each ($P_7, C_5$, kite, paraglider)-free graph $G$ can be polynomially colored with at most $\omega(G)+1$ colors.
\end{corollary}

This also conditionally answers the problem of of Schiermeyer \cite{SR2023} which asks for a constant $c$ such that $\chi(G)\le \omega(G)+c$ for every $(P_7$, diamond)-free graph.

We will prove Theorems~\ref{kite, paraglider} in Section~\ref{2}.

\section{Proof of Theorem~\ref{kite, paraglider}}\label{2}

We begin from the structural properties of $(P_7, C_5$, kite, paraglider)-free graphs that contain a 7-hole.

Let $G$ be a connected imperfect $(P_7, C_5,$ kite, paraglider)-free graph. Since $G$ is imperfect, and an odd antihole on more than five vertices contains paraglider, we have that $G$ contains a 7-hole by Theorem~\ref{perfect}, say $v_1v_2v_3v_4v_5v_6v_7v_1$. Let $L=\{v_1,v_2,\cdots,v_7\}$, and let $R=V(G)\setminus (L\cup N(L))$. For each $i\in \{1, 2, \ldots, 7\}$, we define
\begin{eqnarray*}	
	A_i&=&\{x\in N(L)|~N_L(x)=\{v_i, v_{i+2}\}\},\\
	B_i&=&\{x\in N(L)|~N_L(x)=\{v_i, v_{i+2}, v_{i+3}\}\},\\
	\bar{B}_i&=&\{x\in N(L)|~N_L(x)=\{v_i, v_{i+1}, v_{i+3}\}\},\\
	D_i&=&\{x\in N(L)|~N_L(x)=\{v_i, v_{i+2}, v_{i+3}, v_{i+5} \}\}.
\end{eqnarray*}
Let $I=\{x\in N(L)|~N_L(x)=L\}$,
and let $A=\cup_{i=1}^7 A_i$, $B=\cup_{i=1}^7 B_i$, $\bar{B}=\cup_{i=1}^7 \bar{B}_i$ and $D=\cup_{i=1}^7 D_i$, where the subscript is understood to be modulo 7.
\begin{lemma}\label{N(L)}
$V(G)=L\cup A\cup B\cup \bar{B} \cup D \cup I \cup R$.
\end{lemma}
\pf We need only to verify that $N(L)\subseteq A\cup B\cup \bar{B} \cup D\cup I$. Let $x\in N(L)$, and suppose, without loss of generality, that $x\sim v_1$. It is certain that $|N(x)\cap L|\geq 2$, as otherwise $xv_1v_2v_3v_4v_5v_6$ is an induced $P_7$ of $G$.

Suppose $x\sim v_2$. To avoid an induced $P_7=xv_2v_3v_4v_5v_6v_7$, $|N_L(x)|\ge 3$. If $x\sim v_3$, then $x\sim v_4$ and $x\sim v_7$ to avoid a kite $G[\{x,v_1,v_2,v_3,v_4\}]$ or $G[\{x,v_1,v_2,v_3,v_7\}]$, and consequently $x$ must be complete to $\{v_5,v_6\}$ to avoid a kite $G[\{x,v_2,v_3,v_4,v_5\}]$ or $G[\{x,v_1,v_2,v_6,v_7\}]$, and hence $x\in I$. A similar situation happens if $x\sim v_7$. So, we suppose that $x\not\sim v_3$ and $x\not\sim v_7$. If $x\sim v_5$, then $x\sim v_4$ and $x\sim v_6$ to avoid a 5-hole $xv_2v_3v_4v_5x$ or $xv_5v_6v_7v_1x$, which produce a kite $G[\{x,v_3,v_4,v_5,v_6\}]$. Hence, $x\not\sim v_5$, and so $N_L(x)\subseteq  \{v_1,v_2,v_4,v_6\}$. If $x\not\sim v_4$ then $x\in B_6$. If $x\not\sim v_6$ then $x\in \bar{B}_1$.   Otherwise, $x\in D_6$.

By symmetry, we have that $x\in B\cup \bar{B}\cup D\cup I$ if $x\sim v_2$ or $x\sim v_7$.  Therefore, we suppose that $x\not\sim v_2$ and $x\not\sim v_7$, i.e., $N_L(x)\subseteq \{v_1, v_3, v_4, v_5, v_6\}$

Note that $|N_L(x)|\ge 2$. If $x\not\sim v_3$ and $x\not\sim v_6$, then $G$ has either a 5-hole $xv_1v_2v_3v_4x$ whenever $x\sim v_4$, or a 5-hole $xv_5v_6v_7v_1x$ whenever $x\sim v_5$. This shows that $x\sim v_3$ or $x\sim v_6$. By symmetry, we suppose $x\sim v_3$.

Since $N_L(x)=\{v_1,v_3\}$ implies that $x\in A_1$, we may suppose that $|N_L(x)|\geq 3$. If $x\sim v_5$, then $x\sim v_6$ to avoid a 5-hole  $xv_1v_7v_6v_5x$, and $x\not\sim v_4$ to avoid a kite $G[\{x,v_2,v_3,v_4,v_5\}]$, and thus $x\in D_3$. If $x\not\sim v_5$ and $x\sim v_6$, then $x\sim v_4$ to avoid a 5-hole $xv_3v_4v_5v_6x$, and thus $x\in D_1$. If $x\not\sim v_5$ and $x\not\sim v_6$, then $N_L(x)=\{v_1,v_3,v_4\}$ and so $x\in B_1$. This proves Lemma~\ref{N(L)}. \qed

Next, we prove (M1)$\sim$ (M17), that are some useful properties of $N(L)$. Let $i\in \{1,2,\ldots, 7\}$.

\noindent{\bf (M1)} {\em $A_i$ is a stable set, and no vertex of $A$ has neighbors in $R$}.

Since for any edge $xx'$ of $G[A_i]$, $G[\{x,x',v_i,v_{i+1},v_{i+2}\}]$ is a paraglider, we have that $A_i$ is a stable set.  For any vertex $a\in A_i$, we have an induced $P_6=av_{i+2}v_{i+3}v_{i+4}v_{i+5}v_{i+6}$. It follows from the $P_7$-freeness of $G$ that $N_R(a)=\emptyset$. So, (M1) holds.

\noindent{\bf (M2)} {\em $\max\{|B_i|, |\bar{B}_i|, |D_i|\}\leq 1$}.

Suppose not, let $x$ and $x'$ be two vertices in $B_i$ or $D_i$. If $x\sim x'$ then $G[\{x,x',v_i,v_{i+1},v_{i+2}\}]$ is a paraglider. If $x\not\sim x'$ then $G[\{x,x',v_i,v_{i+2},v_{i+3}\}]$ is a paraglider. Now let $y, y'\in \bar{B}_i$. If $y\sim y'$ then $G[\{y,y',v_{i+1},v_{i+2},v_{i+3}\}]$ is a paraglider. If $y\not\sim y'$, then $G[\{y,y',v_i,v_{i+1},v_{i+3}\}]$ is a paraglider. Therefore, (M2) holds.

{\bf For simplicity, we always take $i=1$ in the following proofs of (M3) $\sim$ (M17)}.

\noindent{\bf (M3)} {\em $A_i$ is complete to $A_{i+1}\cup A_{i-1}$, and anticomplete to $A_{i+3}\cup A_{i-3}$}.

If $x\not\sim y$ for some $x\in A_1$ and $y\in A_2$, then $G$ has an induced $P_7=xv_1v_2yv_4v_5v_6$.
If $x\sim y$ for some $x\in A_1$ and $y\in A_4$, then $xyv_6v_7v_1x$ is a 5-hole. By symmetry, (M3) holds.

\noindent{\bf (M4)} {\em $E(A_i, A_{i+2})$ is a matching}.

Suppose not. Let $e_1=xy_1$ and $e_2=xy_2$ be two edges between $A_1$ and $A_{3}$, where $x\in A_1$ and $y_1$, $y_2\in A_{3}$. Since $y_1\not\sim y_2$ by (M1), we have  a paraglider $G[\{x,y_1,y_2,v_{3},v_{5}\}]$. This shows that (M4) holds.

\noindent{\bf (M5)} {\em Either  $B_i=\emptyset$ or $B_{i+1}\cup B_{i+3}\cup B_{i+4}\cup B_{i+6}=\emptyset$}.

Suppose to its contrary, let $x\in B_1$ and $y\in B_2\cup B_4\cup B_5\cup B_7$. If $y\in B_2$, then $xv_1v_2yv_4x$ is a 5-hole if $x\not\sim y$, and $G[\{x,y,v_2,v_3,v_4\}]$ is a paraglider if $x\sim y$. If $y\in B_4$, then $xv_4yv_7v_1x$ is a 5-hole if $x\not\sim y$, and $G[\{x,y,v_2,v_3,v_4\}]$ is a kite if $x\sim y$. If $y\in B_5$, then $xv_4v_5yv_1x$ is a 5-hole if $x\not\sim y$, and $G[\{x,y,v_1,v_6,v_7\}]$ is a kite if $x\sim y$. If $y\in B_7$, then $xv_3yv_7v_1x$ is a 5-hole if $x\not\sim y$, and $G[\{x,y,v_3,v_4,v_5\}]$ is a kite if $x\sim y$. Therefore, (M5) holds.

\noindent{\bf (M6)} {\em $E(A_i, B_i\cup B_{i+2}\cup B_{i+4}\cup B_{i+5})=\emptyset$, and $E(B_i, B_{i+2}\cup B_{i-2})=\emptyset$}.

First we suppose $x\sim y$ for some $x\in A_1$ and $y\in B_1\cup B_3\cup B_5\cup B_6$. If $y\in B_1\cup B_6$ then $G[\{x,y,v_1,v_2,v_3\}]$ is a paraglider. If $y\in B_3$ then $xyv_6v_7v_1x$ is a 5-hole. If $y\in B_5$ then $G[\{x,y,v_1,v_6,v_7\}]$ is a kite. All are contradictions. Next, suppose that $x\sim y$ for some $x\in B_1$ and $y\in B_3$. Then $G[\{x,y,v_3,v_4,v_5\}]$ is a paraglider. Therefore, (M6) holds.

\noindent{\bf (M7)} {\em Either $D_i=\emptyset$ or $D_{i+1}\cup D_{i+2}\cup D_{i+3}\cup D_{i+4}\cup D_{i+5}\cup D_{i+6}=\emptyset$}.

Suppose to its contrary, we choose that $x\in D_1$ and $y\in D_2\cup D_3\cup D_4$. If $y\in D_2\cup D_4$, then $xv_1v_2yv_4x$ is a 5-hole if $x\not\sim y$, and $G[\{x,y,v_2,v_3,v_4\}]$ is a paraglider if $x\sim y$. If $y\in D_3$, then $xv_4v_5yv_1x$ is a 5-hole if $x\not\sim y$, and $G[\{x,y,v_3,v_4,v_5\}]$ is a paraglider if $x\sim y$. This proves (M7).

\noindent{\bf (M8)} {\em $A_i$ is complete to $B_{i+1}$, and either $A_i=\emptyset$ or $B_{i-1}=\emptyset$}.

First suppose that $x\not\sim y$ for some $x\in A_1$ and $y\in B_2$. Then, $G$ has an induced $P_7=yv_4v_3xv_1v_7v_6$, a contradiction. Next, we suppose that $A_1$ has a vertex $x$ and $B_7$ has a vertex $y$. Then, $G[\{x,y,v_1,v_2,v_3\}]$ is a paraglider if $x\sim y$, and $xv_3yv_7v_1x$ is a 5-hole if $x\not\sim y$. Both are contradictions. Therefore, (M8) holds.

\noindent{\bf (M9)} {\em $E(A_i, D_i\cup D_{i+2}\cup D_{i+4})=\emptyset$, and either $A_i=\emptyset$ or $D_{i+5}\cup D_{i+6}=\emptyset$}.

First suppose that $x\sim y$ for some $x\in A_1$ and $y\in D_1\cup D_3\cup D_5$. Then, $G[\{x,y,v_1,v_2,v_3\}]$ is a paraglider.
Next, we suppose that $A_1$ has a vertex $x$ and $D_6\cup D_7$ has a vertex $y$. If $y\in D_6$, then $xv_3v_4yv_1x$ is a 5-hole if $x\not\sim y$, and $G[\{x,y,v_1,v_2,v_3\}]$ is a paraglider if $x\sim y$. If $y\in D_7$, then $xv_3yv_7v_1x$ is a 5-hole if $x\not\sim y$, and $G[\{x,y,v_1,v_2,v_3\}]$ is a paraglider if $x\sim y$. All are contradictions. This proves (M9).

\noindent{\bf (M10)} {\em $E(B_i, D_{i+5})=\emptyset$, and either $B_i=\emptyset$ or $D\setminus D_{i+5}=\emptyset$}.

If $x\sim y$ for some $x\in B_1$ and $y\in D_6$, then $G[\{x,y,v_1,v_2,v_3\}]$ is a paraglider. Suppose that $B_1$ has a vertex $x$ and $D\setminus D_6$ has a vertex $y$. If $y\in D_1$, then $G$ contains a paraglider $G[\{x,y,v_1,v_3,v_4\}]$ (if $x\not\sim y$) or $G[\{x,y, v_1,v_2,v_3\}]$ (if $x\sim y$).  If $y\in D_2\cup
D_4$, then $xv_1v_2yv_4x$ is a 5-hole if $x\not\sim y$ and $G[\{x,y,v_2,v_3,v_4\}]$ is a paraglider if $x\sim y$. If $y\in D_3\cup D_5$, then $xv_4v_5yv_1x$ is a 5-hole if $x\not\sim y$ and $G[\{x,y, v_1,v_2,v_3\}]$ is a paraglider if $x\sim y$. If $y\in D_7$, then $xv_1v_7yv_3x$ is a 5-hole if $x\not\sim y$ and $G[\{x,y,v_3,v_4,v_5\}]$ is a paraglider if $x\sim y$. All are contradictions. Therefore, (M10) holds.

\noindent{\bf (M11)} {\em $E(\bar{B}_i, \bar{B}_{i+2}\cup \bar{B}_{i+5})=\emptyset$, and either $\bar{B}_i=\emptyset$ or $\bar{B}_{i+1}\cup \bar{B}_{i+3}\cup \bar{B}_{i+4}\cup \bar{B}_{i+6}=\emptyset$}.

If $x\sim y$ for some $x\in \bar{B}_1$ and $y\in \bar{B}_3$, then $G[\{x,y,v_2,v_3,v_4\}]$ is a paraglider. Suppose that $\bar{B}_1$ has a vertex $x$ and $\bar{B}_2\cup \bar{B}_4$ has a vertex $y$. If $y\in \bar{B}_2$, then $xv_4v_5yv_2x$ is a 5-hole if $x\not\sim y$ and $G[\{x,y,v_2,v_3,v_4\}]$ is a paraglider if $x\sim y$. If $y\in \bar{B}_4$, then $xv_4yv_7v_1x$ is a 5-hole if $x\not\sim y$ and $G[\{x,y,v_2,v_4,v_5\}]$ is a kite if $x\sim y$. This proves (M11).

\noindent{\bf (M12)} {\em $A_i$ is complete to $\bar{B}_{i+5}$, and $E(A_i, \bar{B}_{i+1}\cup \bar{B}_{i+2}\cup \bar{B}_{i+4}\cup \bar{B}_{i+6})=\emptyset$}.

If $x\not\sim y$ for some $x\in A_1$ and $y\in \bar{B}_{6}$, then $G$ has an induced $P_7=v_1xv_{3}v_{4}v_{5}v_{6}y$, a contradiction.
Suppose that $x\sim y$ for some $x\in A_1$ and $y\in \bar{B}_2\cup \bar{B}_3\cup \bar{B}_5\cup \bar{B}_7$. Then, $G$ has a paraglider $G[\{x,y,v_1,v_2,v_3\}]$ if $y\in \bar{B}_2\cup \bar{B}_7$, has a kite $G[\{x,y,v_1,v_3,v_4\}]$ if $y\in \bar{B}_3$, and has a 5-hole $xv_3v_4v_5yx$ if $y\in \bar{B}_5$. This proves (M12).

\noindent{\bf (M13)} {\em $E(\bar{B}_i, D_i)=\emptyset$, and either $A_i=\emptyset$ or $\bar{B}_i=\emptyset$}.

If $x\sim y$ for some $x\in \bar{B}_1$ and $y\in D_1$, then $G[\{x,y,v_1,v_{4},v_{5}\}]$ is a kite, a contradiction.
Suppose that $A_1$ has a vertex $x$ and $\bar{B}_1$ has a vertex $y$. Then, $G[\{x,y,v_1,v_2,v_3\}]$ is a paraglider if $x\sim y$, and $xv_3v_4yv_1x$ is a 5-hole if $x\not\sim y$. Therefore, (M13) holds.

\noindent{\bf (M14)} {\em Either $\bar{B}_i=\emptyset$ or $D\setminus D_{i}=\emptyset$}.

Suppose that (M14) is not true. We choose by symmetry that $x\in \bar{B}_1$ and $y\in D\setminus D_1$. If $y\in D_2\cup D_4$, then $G[\{x,y,v_2,v_3,v_4\}]$ is a paraglider if $x\sim y$ and $xv_4yv_7v_1x$ is a 5-hole if $x\not\sim y$. If $y\in D_3\cup D_5$, then $G[\{x,y,v_1,v_2,v_3\}]$ is a paraglider if $x\sim y$ and $xv_4v_5yv_1x$ is a 5-hole if $x\not\sim y$. If $y\in D_6$, then $G[\{x,y,v_2,v_3,v_4\}]$ is a paraglider if $x\sim y$ and $G[\{x,y,v_1,v_2,v_4\}]$ is a paraglider if $x\not\sim y$. If $y\in D_7$, then $G[\{x,y,v_2,v_3,v_4\}]$ is a paraglider if $x\sim y$ and $xv_4v_5yv_2x$ is a 5-hole if $x\not\sim y$. All are contradictions. Therefore, (M14) holds.

\noindent{\bf (M15)} {\em $E(\bar{B}_i, B_i\cup B_{i+1}\cup B_{i+2}\cup B_{i+6})=\emptyset$}.

If it is not true, we choose an edge $xy$ with $x\in \bar{B}_1$ and $y\in B_1\cup B_2\cup B_3\cup B_7$, then $G[\{x,y,v_1,v_4,v_5\}]$ is a kite if $y\in B_1$, $G[\{x,y,v_2,v_3,v_4\}]$ is a paraglider if $y\in B_2\cup B_7$, and $xyv_6v_7v_1x$ is a 5-hole if $y\in B_3$. This proves (M15).

\noindent{\bf (M16)} {\em $\bar{B}_i$ is complete to $B_{i+3}$, and either $\bar{B}_i=\emptyset$ or $B_{i+4}\cup B_{i+5}=\emptyset$}.

First suppose that $x\not\sim y$ for some $x\in \bar{B}_1$ and $y\in B_4$. Then, $v_1xv_4yv_7v_1$ is a 5-hole.
Next, suppose that $\bar{B}_1$ has a vertex $x$ and $B_5\cup B_6$ has a vertex $y$. If $y\in B_5$, then $G[\{x,y,v_1,v_2,v_3\}]$ is a kite if $x\sim y$ and $xv_4v_5yv_1x$ is a 5-hole if $x\not\sim y$. If $y\in B_6$, then $xv_4v_5v_6yx$ is a 5-hole if $x\sim y$ and $G[\{x,y,v_1,v_2,v_4\}]$ is a kite if $x\not\sim y$. All are contradictions. Therefore, (M16) holds.

\noindent{\bf (M17)} {\em If $N(v_i)\subseteq A\cup L$ and $G$ does not have a comparable pair, then $A_{i-1}=\emptyset$}.

Suppose that $N(v_i)\subseteq A\cup L$ and $A_{i-1}\ne\emptyset$. Let $a_{i-1}\in A_{i-1}$. By the definition of $A_{i-1}$, we have that $N(v_i)\subseteq \{v_{i-1},v_{i+1}\}\cup A_{i-2}\cup A_{i}$, $a_{i-1}$ is complete to $\{v_{i-1},v_{i+1}\}$, and $a_{i-1}\not\sim v_i$. Since $a_{i-1}$ is complete to $A_{i-2}\cup A_{i}$ by (M3), we have that $N(v_i)\subseteq N(a_{i-1})$. This proves (M17).

%%%%%%%%%%%%%%%%%%%%%%%%%%%%%%%%%%%%%%%%%%%%%%%%%%%%%%%%%%%%%%%%%%%%%%%%%%%%%%%%%
%%%%%%%%%%%%%%%%%%%%%%%%%%%%%%%%%%%%%%%%%%%%%%%%%%%%%%%%%%%%%%%%%%%%%%%%%%%%20240727

Now, we suppose $G$ is a connected imperfect ($P_7, C_5$, paraglider, kite)-free graph with no clique cutsets and no comparable pairs, and $\delta(G)\geq\omega(G)+1$. By (M2), we always set $B_i=\{b_i\}$ if $B_i\ne\emptyset$, $\bar{B}_i=\{\bar{b}_i\}$ if $\bar{B}_i\ne\emptyset$, and $D_i=\{d_i\}$ if $D_i\ne\emptyset$.

\begin{claim}\label{I}
	If $I\ne\emptyset$, then it is a universal clique.
\end{claim}
\pf Suppose that $I\ne\emptyset$. Let $u\in I$, and let $v\in A\cup
B\cup \bar{B}\cup D$. If $u\not\sim v$, then $G[\{u,v,v_{i+2},v_{i+3},v_{i+4}\}]$ is a kite whenever $v\in A_i$, $G[\{u,v,v_{i+3},v_{i+4},v_{i+5}\}]$ is a kite whenever $v\in B_i\cup \bar{B}_i$, and $G[\{u,v,v_i,v_{i+2},v_{i+3}\}]$ is a paraglider whenever $v\in D_i$. All are contradictions. Therefore, $u\sim v$, and thus $I$ is complete to $A\cup B\cup \bar{B}\cup D$.

Let $G'=G-I$. For integer $h\ge 1$, let $N_h=\{x\in V(G')\setminus L |~d(x, L)=h\}$, where $d(x, L)=\min\{d_{G}(x,x') |~x'\in L\}$. We have that $N_1=A\cup B\cup \bar{B}\cup D$ and thus $I$ is complete to $N_1$. Let $N=\cup_{i\geq1}N_i$, $X=L\cup N$, and $Y=(V(G')\setminus X)\cup I$. Clearly, $V(G)=X\cup Y$.

Obviously, $E(X, V(G')\setminus X)=\emptyset$. We will prove that
\begin{equation}\label{homogeneous}
	\mbox{$X$ is complete to $I$.}
\end{equation}
Suppose not, let $x\in N_h$ for some $h\geq 2$ and $y\in I$ such that $x\not\sim y$. By the definition of $N_h$, there exists an induced path $P=t_ht_{h-1}t_{h-2}\cdots t_1t_0$, where $t_j\in N_j$ for $1\leq j\leq h$, $t_h=x$ and $t_0\in L$. Let $k$ be the smallest integer in $\{2, \ldots, h\}$ such that $t_k\not\sim y$. Since $k\geq 3$ implies that $G[\{y,t_k,t_{k-1},t_{k-2},t_{k-3}\}]$ is a kite, we have that $k=2$, i.e., $y\sim t_0$ and  $y\sim t_1$ but $y\not\sim t_2$. Notice that $N_1=A\cup B\cup \bar{B}\cup D$. We may choose by symmetry that $t_1\in A_1\cup B_1\cup \bar{B}_1 \cup D_1$. Then, either $G[\{v_1,v_2,t_1,t_2,y\}]$ is a kite if $t_1\in A_1\cup B_1\cup D_1$, or $G[\{v_2,v_3,t_1,t_2,y\}]$ is a kite if $t_1\in \bar{B}_1$. This proves~(\ref{homogeneous}).

Now, we have that $X$ is complete to $I$ and anticomplete to $V(G')\setminus X$. If $I$ is not a clique, choose $a$ and $a'$ in $I$ with $a\not\sim a'$, then $G[\{a,a',v_1,v_2,v_4\}]$ is a paraglider. Therefore, $I$ is a clique of $G$. Since $G$ has no clique cutsets, we have that $V(G')=X$, and thus $I$ is a universal clique of $G$. This proves Claim~\ref{I}. \qed

\medskip

From now on, we always suppose that $I=\emptyset$. Then, by Lemma~\ref{N(L)}, $V(G)=A\cup B\cup \bar{B}\cup D\cup L\cup R$. If $B\cup \bar{B}\cup D=\emptyset$, then $N(v_i)\subseteq A\cup L$ for each $i$ and thus $A=\emptyset$ by (M17), i.e., $V(G)=L$, a contradiction to $\delta(G)\ge \omega(G)+1$. Therefore,
\begin{equation}\label{eqa-B-C-D-nonempty-000}
	B\cup \bar{B}\cup D\ne\emptyset, \mbox{ and } \omega(G)\geq3.
\end{equation}

We will show that $\bar{B}\ne\emptyset$ and $D\ne\emptyset$ in the following Claims~\ref{D is not empty} and \ref{C_1 is not empty*}, and then deduce that $G$ must be isomorphic to $F$.
\begin{claim}\label{D is not empty}
	$\bar{B}\ne\emptyset$.
\end{claim}
\pf Suppose to its contrary that $\bar{B}=\emptyset$. By (\ref{eqa-B-C-D-nonempty-000}), $V(G)=A\cup B\cup D\cup L\cup R$, and $B\cup D\ne\emptyset$.

Suppose $D=\emptyset$.  We have then $V(G)=A\cup B\cup L\cup R$ and $B\ne\emptyset$. Without loss of generality, we suppose that $B_1=\{b_1\}$. By (M5) and (M8), $A_2\cup B_2\cup B_4\cup B_5\cup B_7=\emptyset$, and thus $V(G)=(\cup_{i\ne 2} A_i)\cup \{b_1\}\cup B_3\cup B_6\cup L\cup R$. If $B_6=\emptyset$, then $b_1$ is complete to $\{v_1,v_3\}\cup A_7$ by (M8), which forces a comparable pair $(b_1, v_2)$ as $N(v_2)=\{v_1,v_3\}\cup A_7\subseteq N(b_1)$. Therefore, $B_6\ne\emptyset$, which implies that $B_3=\emptyset$ by (M5), and $A_7=\emptyset$ by (M8). Consequently, $N(v_5)\subseteq L\cup A$ and $N(v_7)\subseteq L\cup A$. By (M17), we have that $A_4\cup A_6=\emptyset$, and hence $V(G)=A_1\cup A_3\cup A_5\cup \{b_1, b_6\}\cup L\cup R$. Now, $N(v_7)=\{v_1,v_6\}\cup A_5$. Since $b_6$ is complete to $\{v_1,v_6\}\cup A_5$ by (M8), and since $b_6\not\sim v_7$, we get a comparable pair $(b_6, v_7)$. So, $D\ne\emptyset$.

Without loss of generality, we suppose that $D_1\ne\emptyset$. By (M7), (M9) and (M10), we have that $A_2\cup A_3\cup (B\setminus B_3)\cup(D\setminus D_1)=\emptyset$, and thus $V(G)=(\cup_{i\ne 2, 3} A_i)\cup B_3\cup \{d_1\}\cup L\cup R$. Since $N(L)\setminus A\subseteq\{b_3, d_1\}$, we have that $N(v_j)\subseteq A\cup L$ for each $j\in\{2,7\}$, and thus $A_1\cup A_6=\emptyset$ by (M17). Now, we have $V(G)=A_4\cup A_5\cup A_7\cup B_3\cup \{d_1\}\cup L\cup R$, which implies $d(v_1)=3\leq \omega(G)$, a contradiction. This proves that $\bar{B}\ne\emptyset$, and thus Claim~\ref{D is not empty} holds. \qed

Now, we suppose, without loss of generality, that $\bar{B}_1=\{\bar{b}_1\}$. By (M11), (M13), (M14) and (M16), we have that $A_1\cup B_5\cup B_6\cup (\bar{B}\setminus(\bar{B}_1\cup \bar{B}_3\cup \bar{B}_6))\cup (D\setminus D_1)=\emptyset$. Hence,
\begin{equation}\label{eqa-A1-emptyset}
	\mbox{$V(G)=(\cup_{2\le i\le 7} A_i)\cup(\cup_{i\not\in \{5, 6\}} B_i)\cup \{\bar{b}_1\}\cup \bar{B}_3\cup \bar{B}_6\cup  D_1\cup L\cup R$.}
\end{equation}
\begin{claim}\label{C_1 is not empty*}
	$D=D_1\ne\emptyset$.
\end{claim}
\pf Suppose to its contrary that $D_1=\emptyset$. Then, $V(G)=(\cup_{2\le i\le 7} A_i)\cup(\cup_{i\not\in \{5, 6\}} B_i)\cup \{\bar{b}_1\}\cup \bar{B}_3\cup \bar{B}_6\cup L\cup R$.  By (M11), either $\bar{B}_3=\emptyset$ or $\bar{B}_6=\emptyset$, and so $\bar{B}\in \{\{\bar{b}_1\}, \{\bar{b}_1, \bar{b}_3\}, \{\bar{b}_1,\bar{b}_6\}\}\}$.  We  show first that
\begin{equation}\label{D_1, D_3}
	\bar{B}\ne \{\bar{b}_1,\bar{b}_3\}.
\end{equation}
Suppose to its contrary that $\bar{B}=\{\bar{b}_1,\bar{b}_3\}$. By (M11), (M13) and (M16), $A_3\cup B_1\cup B_7=\emptyset$ and $\bar{b}_1\not\sim \bar{b}_3$, and so $V(G)=(\cup_{i\ne 1,3} A_i)\cup(\cup_{i\in \{2,3,4\}} B_i)\cup \{\bar{b}_1,\bar{b}_3\}\cup L\cup R$. By (M5), we have that one of $B_2$ and $B_3$ is empty, and one of $B_3$ and $B_4$ is empty.

Suppose $B_3\ne\emptyset$. Then, $B_2\cup B_4=\emptyset$, and so $B=B_3=\{b_3\}$. Notice that $A_4=\emptyset$ by (M8). Since $N(L)\setminus A=\{b_3,\bar{b}_1,\bar{b}_3\}$, we have that  $N(v_7)\subseteq A\cup L$, and thus $A_6=\emptyset$ by (M17). Now, $A=A_2\cup A_5\cup A_7$, which leads to a contradiction that $d(v_1)=3\leq \omega(G)$.

Therefore, $B_3=\emptyset$, and so $B=B_2\cup B_4$.

If $B=B_2\cup B_4=\{b_2,b_4\}$, then $A_3\cup A_5=\emptyset$ by (M8), which leads to a contradiction that $d(v_3)=d(v_5)=3\leq \omega(G)$ as $A=A_2\cup A_4\cup A_6\cup A_7$.

If $B=\emptyset$, we have that $N(L)\setminus A=\{\bar{b}_1, \bar{b}_3\}$ and thus $N(v_j)\subseteq A\cup L$ for each $i\in\{5,7\}$, then $A_4\cup A_6=\emptyset$ by (M17), which leads to a contradiction that $d(v_1)=d(v_3)=d(v_6)=3\leq \omega(G)$ as $A=A_2\cup A_5\cup A_7$.

Suppose that  $B=B_2=\{b_2\}$. Since $N(L)\setminus A=\{b_2,\bar{b}_1,\bar{b}_3\}$, we have that $N(v_7)\subseteq A\cup L$ and thus $A_6=\emptyset$ by (M17), which leads to a contradiction that $d(v_1)=d(v_3)=3\leq \omega(G)$ as $A=A_2\cup A_4\cup A_5\cup A_7$.

Suppose that $B=B_4=\{b_4\}$. Then, $A_5=\emptyset$ by (M8). Since $N(L)\setminus A=\{b_4,\bar{b}_1,\bar{b}_3\}$, we have that $N(v_5)\subseteq A\cup L$ and thus $A_4=\emptyset$ by (M17), which leads to a contradiction that $d(v_5)=2<\omega(G)$ as $A=A_2\cup A_6\cup A_7$.
This completes the proof of (\ref{D_1, D_3}).

Using a similar argument, we can show that  $\bar{B}\ne \{\bar{b}_1, \bar{b}_6\}$.

Therefore, $\bar{B}=\{\bar{b}_1\}$. Recall that $V(G)=(\cup_{2\le i\le 7} A_i)\cup(\cup_{i\not\in \{5, 6\}} B_i)\cup \{\bar{b}_1\}\cup L\cup R$ and $B=\cup_{i\not\in\{5, 6\}} B_i$. We show next that
\begin{equation}\label{eqa-B4-B7-000}
	B=B_4\cup B_7.
\end{equation}
Suppose $B_1\ne\emptyset$. Then, $A_2\cup B_2\cup B_4\cup B_5\cup B_7=\emptyset$ by (M5) and (M8), and thus $B=B_1\cup B_3$. Since $N(L)\setminus A\subseteq\{b_1,b_3,\bar{b}_1\}$, we have that $N(v_7)\subseteq A\cup L$, and thus $A_6=\emptyset$ by (M17). If $B_3\ne\emptyset$, then $A_4=\emptyset$ by (M8). If $B_3=\emptyset$, then $N(v_5)\subseteq A\cup L$ and so $A_4=\emptyset$ by (M17). Both lead to a contradiction that $d(v_6)\leq3\leq\omega(G)$ as $A_4\cup A_6=\emptyset$. This shows that $B_1=\emptyset$, and hence $B=B_2\cup B_3\cup B_4\cup B_7$.

If $B_2\ne\emptyset$, then $A_3\cup B_3=\emptyset$ by (M5) and (M8), and thus $A=A\setminus(A_1\cup A_3)$ and $B=B_2\cup B_4\cup B_7$, which leads to a contradiction that $d(v_3)=3\leq \omega(G)$. So, we have further that $B_2=\emptyset$, and $B=B_3\cup B_4\cup B_7$.

Suppose $B_3\ne\emptyset$. Then, $A_4\cup B_4\cup B_7=\emptyset$ by (M5) and (M8), and thus $B=B_3=\{b_3\}$. Since $N(L)\setminus A=\{b_3,\bar{b}_1\}$, we have that $N(v_7)\subseteq A\cup L$ and thus $A_6=\emptyset$ by (M17). Now, $A=A_2\cup A_3\cup A_5\cup A_7$, which leads to a contradiction that $d(v_1)=d(v_6)=3\leq \omega(G)$. This shows that $B_3=\emptyset$,  and completes the proof of (\ref{eqa-B4-B7-000}).

Recall that $V(G)=(\cup_{2\le i\le 7} A_i)\cup B_4\cup B_7\cup \{\bar{b}_1\}\cup L\cup R$. By (M5), we have that one of  $B_4$ and $B_7$ must be empty.

If $B_4\cup B_7=\emptyset$, then $N(L)\setminus A=\{\bar{b}_1\}$, and so $A_2\cup A_4\cup A_5\cup A_6=\emptyset$ by (M17). But $d(v_1)=3\leq\omega(G)$ as $A=A_3\cup A_7$. Therefore, exact one of $B_4$ and $B_7$ is empty.

Suppose $B_4=\{b_4\}$ and $B_7=\emptyset$. Then, $A_5=\emptyset$ by (M8). Since $N(L)\setminus A=\{b_4, \bar{b}_1\}$, we have that $A_2\cup A_4=\emptyset$ by (M17), and so $A=A_3\cup A_6\cup A_7$ and  $N(v_5)=\{v_4, v_6\}\cup A_3$. Since $b_4$ is complete to $\{v_4,v_6\}\cup A_3$ and $b_4\not\sim v_5$ by (M8), we have  a contradiction that $N(v_5)\subseteq N(b_4)$.

Therefore, $B_4=\emptyset$, $B_7=\{b_7\}$.
Since $N(L)\setminus A=\{b_7, \bar{b}_1\}$, we have that $N(v_5)\cup N(v_6)\subseteq A\cup L$ and thus $A_4\cup A_5=\emptyset$ by (M17). Hence, $A=A_2\cup A_3\cup A_6\cup A_7$.

Since $A_6=\emptyset$ leads to $d(v_6)=2<\omega(G)$, let $a_6\in A_6$. By (M3) and (M8), $a_6$ is complete to $\{v_1,v_6\}\cup A_7\cup B_7$ and $a_6\not\sim v_7$; then $N(v_7)=\{v_1,v_6\}\cup A_7\cup B_7\subseteq N(a_6)$. So, Claim~\ref{C_1 is not empty*} holds.\qed

\medskip

By (\ref{eqa-A1-emptyset}) and Claim~\ref{C_1 is not empty*}, we have that $D=\{d_1\}$. Then, $A_2\cup A_3\cup (B\setminus B_3)\cup (\bar{B}\setminus \bar{B}_1)=\emptyset$ by (M9), (M10) and (M14), which leads to
$\mbox{$V(G)=A_4\cup A_5\cup A_6\cup A_7\cup B_3\cup \{\bar{b}_1, d_1\}\cup L\cup R$}.$

If $B_3=\emptyset$, then leads to a contradiction that $d(v_3)=3\leq\omega(G)$. Therefore, $B=B_3=\{b_3\}$.

By (M8), $A_4=\emptyset$. Since $N(L)=A\cup \{b_3,\bar{b}_1,d_1\}$, we have that $N(v_7)\subseteq A\cup L$ and thus $A_6=\emptyset$ by (M17), which leads to \begin{equation}\label{eqa-A5-A7-B3-C--D1}
	\mbox{$V(G)=A_5\cup A_7\cup \{b_3, \bar{b}_1, d_1\}\cup L\cup R$.}
\end{equation}

By (M10), (M13) and (M15), we have that $\{b_3, \bar{b}_1, d_1\}$ is a stable set. Furthermore, $E(b_3, A_5)=E(\bar{b}_1, A_7)=\emptyset$ by (M6) and (M12). We show next that
\begin{equation}\label{cutset}
	\mbox{$A_5\ne \emptyset$, and $A_7\ne\emptyset$. }
\end{equation}

Suppose that $A_5=\emptyset$. If $\bar{b}_1$ has no neighbor in $R$, then $d(\bar{b}_1)=3\leq\omega(G)$, a contradiction. So, we suppose that $\bar{b}_1$ has a neighbor, say $x$, in $R$. Let $Q$ be the component of $R$ which contains $x$. Notice that $x\bar{b}_1v_4v_5v_6v_7$ is an induced $P_6$. To avoid an induced $P_7$, $\bar{b}_1$ must be complete to $V(Q)$. Let $x'$ be an arbitrary vertex of $Q$. By (M1), we have that $E(x', A_5\cup A_7)=\emptyset$. If $x'\sim d_1$, then $x'\bar{b}_1v_2v_3d_1x'$ is a 5-hole. If $x'\sim b_3$, then  $x'\bar{b}_1v_4v_3b_3x'$ is a 5-hole. Therefore, $E(V(Q), N(L)\setminus \{\bar{b}_1\})=\emptyset$, which leads to a contradiction that $\bar{b}_1$ is a cut vertex. With a similar argument by considering the relation between $b_3$ and $R$, we can deduce a same contradiction if $A_7=\emptyset$. This proves (\ref{cutset}).

By (\ref{eqa-A5-A7-B3-C--D1}) and (\ref{cutset}), we have that $V(G)=A_5\cup A_7\cup \{b_3, \bar{b}_1, d_1\}\cup L\cup R$, $A_5\ne\emptyset$, and $A_7\ne \emptyset$.

If $E(d_1, A_5)\ne\emptyset$, choose $x\in N_{A_5}(d_1)$, then $d_1v_3b_3v_5xd_1$ is a 5-hole, a contradiction. A similar contradiction happens if $E(d_1, A_7)\ne\emptyset$. Since $A_5\ne\emptyset$, let $a\in A_5$, then $a\sim \bar{b}_1$ to forbid an induced $P_7=\bar{b}_1v_2v_3d_1v_6v_5a$. Similarly, we can show that $b_3$ is complete to $A_7$. Therefore,
\begin{equation}\label{A_5}
	\mbox{$\bar{b}_1$ is complete to $A_5$, $b_3$ is complete to $A_7$, and $E(d_1, A_5\cup A_7)=\emptyset$.}
\end{equation}

We show next that $R=\emptyset$. Since $G$ is $C_5$-free and $\bar{b}_1v_2v_3d_1$, $\bar{b}_1v_4v_3b_3$ and $d_1v_4v_5b_3$ are all induced $P_4$, we have that $N(\bar{b}_1)\cap N(d_1)\cap R=\emptyset$, $N(\bar{b}_1)\cap N(b_3)\cap R=\emptyset$ and $N(d_1)\cap N(b_3)\cap R=\emptyset$. Moreover, since both $\bar{b}_1v_2v_3d_1v_6v_7$ and $b_3v_5v_4d_1v_1v_7$ are induced $P_6$, following from the $P_7$-freeness of $G$, we have that $N_R(b_3)=N_R(\bar{b}_1)=\emptyset$.

Suppose that $R\ne \emptyset$. Then, $N_R(d_1)\ne\emptyset$ by (M1). Let $r\in N_R(d_1)$ and $t\in A_5\cup A_7$.   If $t\in A_5$, then $t\sim \bar{b}_1$ by (\ref{A_5}) which forces an induced $P_7=rd_1v_3v_2\bar{b}_1tv_7$. If $t\in A_7$, then $t\sim b_3$ by (\ref{A_5}) which forces an induced $P_7=rd_1v_4v_5b_3tv_7$. Therefore, $R=\emptyset$.

Since $A_5\ne\emptyset$ and $A_7\ne\emptyset$, let $x\in A_5$ and $y\in A_7$, then $\bar{b}_1\sim x$ and $b_3\sim y$ by (\ref{A_5}), and so $x\sim y$ to avoid an induced $P_7=d_1v_4\bar{b}_1xv_7yb_3$. Now, $A_5$ is complete to $A_7$. Consequently, $|A_5|=|A_7|=1$ by (M4), and thus $G$ must be isomorphic to $F$. This proves Theorem~\ref{kite, paraglider}. \qed

\renewcommand{\baselinestretch}{1.1}

\noindent{\bf Statements and Declarations} The authors declare that they have no known competing financial interests or personal relationships that could have appeared to influence the work reported in this paper. The second author is supported by NSFC 11931006. All authors contributed to the study conception and design. Material preparation and analysis were performed by Ran Chen and Baogang Xu. The first draft of the manuscript was written by Ran Chen and all authors commented on previous versions of the manuscript. All authors read and approved the final manuscript.

\end{document}